\documentclass[article,12pt]{elsarticle}
\usepackage{xcolor}
\usepackage{algorithm}
\usepackage{algorithmic}
\usepackage{multirow}
\usepackage{caption}
\usepackage{array}
\usepackage{booktabs}
\usepackage{graphicx}
\usepackage{subcaption}
\usepackage{geometry}

\usepackage{amssymb}
\usepackage{amsmath}

\makeatletter
\def\ps@pprintTitle{%
   \let\@oddhead\@empty
   \let\@evenhead\@empty
   \let\@oddfoot\@empty
   \let\@evenfoot\@empty
}
\begin{document}

\begin{frontmatter}



\title{AttenMfg: An Attention Network Based Optimization Framework for Sensor-Driven Operations \& Maintenance in Manufacturing Systems} 


\author{Iman Kazemian$^\dagger$, Murat Yildirim$^\dagger$ $^*$, and Paritosh Ramanan$^\ddagger$\\
			$^\dagger$ Industrial and Systems Engineering, Wayne State University, Detroit, US \\
             $^\ddagger$ Industrial Engineering \& Management, Oklahoma State University, Stillwater, US \\
             $^*$ Corresponding author email: murat@wayne.edu} 


\begin{abstract}
Operations and maintenance (O\&M) scheduling is a critical problem in leased manufacturing systems, with significant implications for operational efficiency, cost optimization, and machine reliability. Solving this problem involves navigating complex trade-offs between machine-level degradation risks, production throughput, and maintenance team logistics across multi-site manufacturing networks. Conventional approaches rely on large-scale Mixed Integer Programming (MIP) models, which, while capable of yielding optimal solutions, suffer from prolonged computational times and scalability limitations. To overcome these challenges, we propose \textit{AttenMfg}, a novel decision-making framework that leverages multi-head attention (MHA), tailored for complex optimization problems. The proposed framework incorporates several key innovations, including constraint-aware masking procedures and novel reward functions that explicitly embed mathematical programming formulations into the MHA structure. The resulting attention-based model (i) reduces solution times from hours to seconds, (ii) ensures feasibility of the generated schedules under operational and logistical constraints, (iii) achieves solution quality on par with exact MIP formulations, and (iv) demonstrates strong generalizability across diverse problem settings. These results highlight the potential of attention-based learning to revolutionize O\&M scheduling in leased manufacturing systems.
\end{abstract}




\end{frontmatter}



\section{Introduction}
\label{intro}
Leased manufacturing systems are production environments where equipment and machinery are rented rather than owned by the manufacturing company \cite{lyu2024product}. This approach allows manufacturers to minimize upfront capital costs and transfer some maintenance responsibilities to the lessors \cite{liu2024optimal}. According to McKinsey \& Company, incorporating leasing into manufacturing operations enables 12\% to 18\% reduction in logistics costs, 20\% to 24\% decrease in quality-related expenses, and 30\% drop in maintenance costs \cite{mckinsey_indirect_2023}. Nonetheless, the success of this approach hinges on efficient asset management and maintenance scheduling, which is typically handled by the lessor. Scheduling maintenance by the lessor is a complex task due to numerous factors that must be taken into consideration, such as resource availability, operational constraints, and cost-effectiveness of different companies that operate the leased machines. To this end, sensor-driven maintenance approaches can offer significant visibility into condition of the machines and potential failure risks, and enable proactive management of downtime, operations and logistics \cite{deloitte_maintenance_2017}.

Despite its advantages, the adoption of sensor-driven maintenance in leased manufacturing systems remains limited due to three key challenges. First, existing sensor-driven maintenance methods predominantly focus on single assets or multiple assets with only limited and static interactions \cite{schutz2013maintenance,taghipour2012optimum}. In reality, these interactions are highly complex and dynamic, influenced by time-varying demand and asset availability. To unlock the benefits of sensor data, we need integrated models that can faithfully represent the implications of different asset-level failure predictions onto their system-level consequences. Second, a lack of model adaptivity hinders implementation in evolving operational environments. Traditional approaches require specialized expertise to develop and update models in response to changing conditions, making adaptation both costly and resource-intensive. This challenge is particularly limiting for small and medium-sized manufacturers that often lack the necessary skilled labor and capital \cite{deloitte_maintenance_2017}. Finally, managing maintenance across multiple locations introduces logistical complexities. Variations in service standards, scheduling constraints, costs, and reliability across different sites make it difficult to implement a standardized and efficient maintenance strategy. Addressing these challenges is essential for unlocking the full potential of sensor-driven maintenance in leased manufacturing systems \cite{si2022service}.

In this study, we propose an attention-based operations and maintenance framework to address these three fundamental challenges in leased manufacturing systems. The proposed model, \textit{AttenMfg: Attention-Based Operation \& Maintenance Model for Manufacturing Systems}, leverages attention mechanisms to enhance decision-making by capturing complex asset interactions, improving model adaptivity, and optimizing maintenance strategies across dynamic operational environments. Our approach integrates key aspects of operations and maintenance planning, including (i) optimization of sensor-driven maintenance schedules for leased machines based on their real-time degradation states, (ii) efficient routing and coordination of maintenance crews across multiple manufacturing sites, and (iii) minimization of costs associated with machine unavailability and disruptions to planned production levels. The model dynamically incorporates these factors through attention mechanisms, enabling it to handle the complexity of large-scale systems and account for operational constraints. The proposed model also explicitly embeds feasibility within the model's structure by employing constraint-aware masking procedures, which ensures practical applicability and robustness. The attention mechanism enhances the model's ability to interpret and prioritize critical data, leading to improved decision-making, computational efficiency, and solution quality. This novel approach underscores the potential of machine learning in advancing operations research methodologies for maintenance and scheduling challenges in manufacturing.

The remainder of the paper is organized as follows: Section~\ref{litrature} contains the literature review of maintenance in leased manufacturing systems and its solution approaches.
Section~\ref{description} provides a detailed problem description, elucidating the specific characteristics and mathematical formulation of the maintenance scheduling problem in leased manufacturing systems. Section~\ref{MIP} provides a mathematical formulation to solve operation and maintenance scheduling in leased manufacturing systems. Section~\ref{MHA} presents the transformation of the maintenance scheduling problem into a learning problem, where we define the key inputs, outputs, and features necessary for training the model. In section~\ref{attention}, we introduce the attention model, along with its architecture, the role of multi-head attention, and how it effectively captures temporal and non-temporal features relevant to machine maintenance and logistics. Finally, section~\ref{experiment} discusses the experiments conducted to evaluate the proposed model's performance.

\section{Literature Review}
\label{litrature}
We provide a foundation for our proposed study by reviewing two key areas. First, we review the literature on modeling maintenance scheduling in leased manufacturing systems, covering traditional methods, opportunistic strategies, models that consider inter-machine dependencies, economic considerations, and the role of sensor-enabled predictive maintenance. Second, we examine algorithms and approaches to solve these problems, surveying operations research (OR), heuristics, and machine learning (ML) techniques. 


Operations and maintenance scheduling in manufacturing systems has a longstanding history, initially focusing on single-machine models (\cite{chang2011joint,chen2009minimizing,schutz2013maintenance,yeh2009optimal,liu2024optimal}). Although these studies offer foundational insights, they may not always represent modern production environments that have a wide range of operational interdependencies. Opportunistic maintenance strategies emerged to address system-wide interdependencies by repairing additional machinery when significant degradation or failure occurs (\cite{bedford2011signal,feng2023multi,taghipour2012optimum,xia2012dynamic,zhou2009opportunistic}), but they often overlook the broader impact of failures on factory output. Some models incorporate economic dependencies (\cite{koochaki2012condition,wang2010effective,xia2015production,yang2008maintenance,zhou2012preventive,liu2025optimal,sharafali2019optimal,polotski2019joint,rivera2016production,rivera2018subcontracting}), yet complexity rises further for equipment lessors who must manage throughput and availability across multiple sites. This literature primarily schedules maintenance based on general failure statistics rather than real-time machine data. As a result, decisions rely on historical averages rather than asset-specific conditions, limiting precision and potentially increasing downtime and operational disruptions.

Recent advances in sensors, data processing, and storage have enabled new methods that utilize real-time data, such as vibration and temperature to enhance asset condition monitoring and failure prediction \cite{jung2017vibration}. However, integrating these sensor-driven insights into operations and maintenance scheduling remains a significant challenge. Existing literature largely assumes static and limited asset interactions \cite{si2022reliability,jia2020joint}, restricting the ability to model complex and time-varying interdependencies. While some studies incorporate maintenance dependencies (\cite{bedford2011signal,camci2015maintenance}), they do not integrate both economic and maintenance interdependencies, which is essential for leased systems. The work in \cite{karakaya2024sensor} addresses this gap by developing a sensor-driven model that optimizes failure rate predictions, maintenance team routing, and production flow across sites. Despite outperforming traditional time-based methods, these frameworks still face challenges in adaptation to diverse operating conditions and in overcoming the computational bottlenecks of mixed integer programming.

After formulating operations and maintenance models, the next challenge is determining how to solve them efficiently. The maintenance scheduling problem is a well-known NP-hard combinatorial optimization challenge \cite{geurtsen2023integrated}, typically with additional complexities introduced by stochastic variables, dynamic costs, and production constraints. The literature addresses these challenges through three main approaches: mixed-integer linear programming (MILP), heuristic/metaheuristic techniques, and machine learning (ML)-based strategies, each with distinct strengths and limitations.

Mixed integer programming (MIP) models, typically solved with commercial tools, can potentially provide optimal solutions; but remain computationally expensive and challenging for large-scale problem settings \cite{chen2020approximate, azadeh2015solving, karakaya2024sensor}. Heuristic methods, such as Genetic Algorithms (GA) \cite{naderi2009scheduling}, Particle Swarm Optimization (PSO) \cite{azadeh2015solving}, and Simulated Annealing (SA) \cite{naderi2009scheduling}, offer faster, near-optimal solutions but often violate problem constraints, result in unpredictable solution quality, and reduce reliability.

In recent years, machine learning techniques have been explored to enhance optimization by approximating lower bounds \cite{baltean2018selecting}, improving branching strategies \cite{alvarez2017machine}, or directly generating decisions for structured problems like the Traveling Salesman Problem (TSP) using pointer networks \cite{vinyals2015pointer} and reinforcement learning \cite{kool2018attention}. However, these methods typically struggle with feasibility, scalability, and generalization in complex maintenance scheduling problems. To address these problems, Attention-based models have gained traction in optimization tasks \cite{vaswani2017attention}, including UAV \cite{wan2024deep} and crane scheduling \cite{jin2023deep} which employed multi-head attention without accounting for the dynamic costs associated with scheduling. MHA transformers excel in sequential decision-making for routing problems \cite{kool2018attention, liu2022good}. \cite{kazemian2024attention} applied MHA to wind energy maintenance, integrating operations research models into a learning-based framework. Their approach improves computational efficiency, feasibility under constraints, and transferability, outperforming traditional MIP methods.

Our contributions are as follows:
\begin{enumerate}
    \item We develop a novel multi-head attention model, \textit{AttenMfg}, designed to optimize sensor-driven operations and maintenance (O\&M) scheduling for leased manufacturing systems. The proposed model effectively integrates (i) real-time sensor-based predictions of machine degradation with (ii) optimal scheduling decisions, enhancing the ability to process high-dimensional, dynamic data and improving scheduling quality.

    \item The feasibility and optimality of the proposed framework are achieved by generalizing the mappings from the mathematical programming formulation’s objective and constraints into the architecture of the multi-head attention model. Key innovations include:
    \begin{itemize}
        \item (i) Constraint-aware masking procedures that ensure a one-to-one mapping between problem-specific operational constraints (e.g., maintenance limits, team routing) and their corresponding masked representations during training and inference.
        \item (ii) Novel reward functions that incorporate maintenance costs, production penalties, and team relocation expenses, enabling the model to holistically evaluate O\&M decisions across large-scale, multi-site manufacturing networks.
    \end{itemize}

    \item We present a comprehensive end-to-end framework that seamlessly integrates attention-based learning and optimization techniques to generate feasible, cost-effective maintenance schedules for leased systems.

    \item We evaluate the generalizability of our proposed model across varying operational scenarios, including both fixed-site and random-site configurations. Our experiments reveal two notable findings: 
    \begin{itemize}
        \item (i) The model achieves superior results when evaluated on datasets with the same site structure as the training datasets.
        \item (ii) Training on larger and more complex datasets enables robust performance when tested on a diverse set of smaller-scale problems, demonstrating strong generalization and transfer learning capabilities.
    \end{itemize}
\end{enumerate}

\section{Problem Description}
\label{description}
We study the fundamental problem of O\&M scheduling in leased manufacturing systems, where the focus is to maximize operational efficiency and profitability while addressing constraints related to production (e.g., maintaining throughput capacity to meet demand) and maintenance (e.g., routing maintenance teams effectively). In this context, we use real-time sensor data to determine optimal schedules for the routing of maintenance teams and the maintenance of machines across the manufacturing network.

The model incorporates two types of maintenance actions: corrective maintenance (CM) and preventive maintenance (PM). Corrective maintenance is performed to restore machines that have already failed, with the machines remaining non-operational until repairs are completed. Preventive maintenance, on the other hand, is a proactive maintenance taken before failures occur, causing downtime only during the maintenance activity itself. The model ensures that these maintenance actions are balanced to minimize disruptions to production while optimizing machine reliability.

A critical aspect of the problem is maintaining sufficient throughput capacity to meet production requirements for a single product. Shutting down a subset of machines, depending on their manufacturing facility's layout, can significantly affect the production capacity of product. To address this, the model integrates machine availability with production demands, ensuring that maintenance actions have minimal adverse effects on throughput.

The dependency between machine maintenance and the routing of maintenance teams is another key consideration. Maintenance teams must be present at a facility before any machines at that location can undergo maintenance. Moreover, the number of machines a team can maintain is limited to a certain number and also we have the cost associated with moving teams and by adjusting the costs of team visits, the model adaptively determines the extent of opportunistic maintenance actions.

Additionally, we have the operational and production impacts of maintenance through specific variables that impose penalties for machine downtime and production interruptions. This penalties enables the lessor, acting as the maintenance scheduler, to dynamically balance maintenance decisions, machine failure risks, operational impacts, and time-varying production demands. The scheduling problem thus considers not only the routing of maintenance teams but also the broader operational outcomes for manufacturers within the network. The optimized maintenance schedules derived from this model are subsequently implemented to support the manufacturing facilities efficiently.

\begin{figure}[h]
    \centering
    \includegraphics[width=15cm]{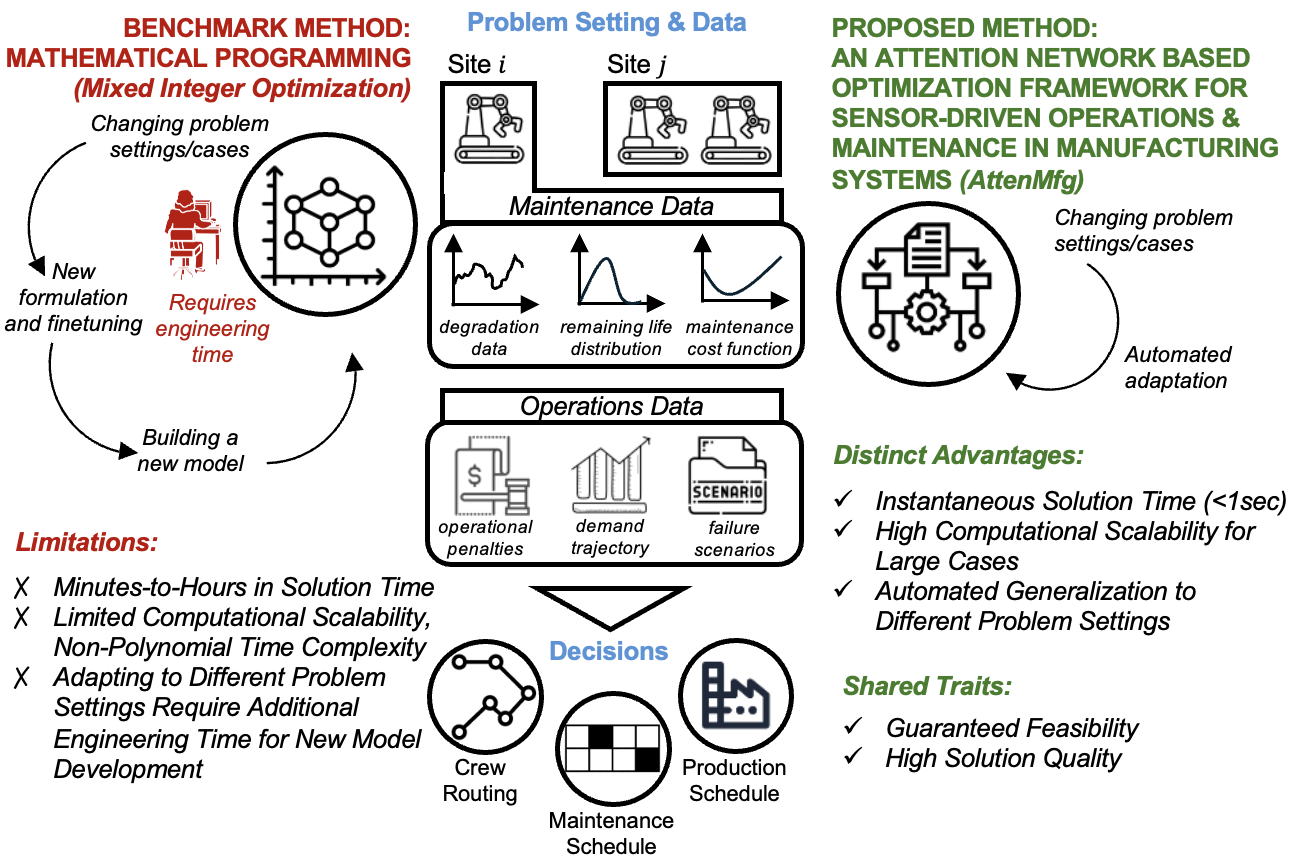}
    \caption{Comparison of the mixed integer optimization model approach and the MHA-based method for operation and maintenance scheduling in leased manufacturing system. The benchmark model struggles with adaptability and scalability, while the MHA approach leverages real-time data for robust, scalable, and efficient maintenance scheduling across diverse scenarios.}
    \label{fig:why}
\end{figure}

As shown in figure~\ref{fig:why}, one of the conventional approaches to solve this problem is the Mixed-Integer Programming (MIP) model, which requires an engineer to construct the model and solve it manually. However, this approach comes with several disadvantages, including the need for specialized expertise, slow response times, and limited scalability, particularly when handling large datasets or adapting to changing maintenance needs.

In this paper, we propose the MHA transformation algorithm to address these challenges. The advantages of this approach include robustness in ensuring consistent and satisfactory performance across different scenarios, quick response times to adapt to dynamic maintenance requirements, high generalizability to perform effectively across diverse industrial settings, and scalability to handle increasing complexity and larger datasets without performance loss.

\section{Lease manufacturing Maintenance Scheduling Problem as a Mathematical Formulation}
\label{MIP}

The conventional approach to addressing the problem described in section~\ref{description} is to formulate the model as a Mixed-Integer Programming (MIP) model. In this section, we develop an MIP model, including the relevant decision variables, parameters, objective function and constraints to serve as a benchmark for our proposed method. This formulation will also help with explaining how the proposed model compares to this benchmark from a modeling perspective.

\subsection{Decision Variables}
\begin{itemize}
     \item \( z_{m,t}\):= 1, if maintenance scheduled for machine \(m\) at period \( t \); otherwise 0.
     \item \( \phi_{l,t}\):= 1, if maintenance team is in manufacturing facility \(l\) at period \( t \); otherwise 0.
     
    \item \( \delta_t \):= 1, if the location of maintenance team changes its location at period \( t \); otherwise 0. This variable will be used to incur a cost for crew transportation.
    
    \item \( d^{u}_{m,t} \): Expectation for the amount of unsatisfied demand at machine \(m\) at period \( t \).
    \item \( \gamma^{u,s}_{m,t} \): Amount of unsatisfied demand at machine \(m\) at period \( t \) in scenario \(s\).
    \item \( \lambda_{m,t}^s \):= Production level for machine \(m\) in period \(t\) in scenario \(s\) due to failure or during the maintenance period.
\end{itemize}

\subsection{Parameters}

\begin{itemize}
     \item \( L \):= Set of manufacturing facilities.
     \item \( M \):= Set of machines.
    \item \( M_n \):= Set of machines in manufacturing facility \(n\).
    \item \( J \):= maximum number of maintenance per period.
    \item \( T \):= Time periods.
    \item \( S \):= Set of scenarios representing possible failure times for the machines.
    \item \( P^f \):= Penalty cost per unit time if a machine stays idle.
    \item \( P^d \):= Penalty cost per unit time per unit product if a machine in manufacturing facility not be able to satisfy the demand.
    \item \( \Delta \):= The cost of changing location for maintenance team.
    \item \( C_{m,t} \):= Estimated maintenance cost for machine \(m\) in period \(t\).
    \item \( C_m^{pf} \):= Penalty cost that lessor should pay due to the failure of machine \(m\) (Unexpected failure cost).
    \item \(P^s_{m,t} \):= Maximum production rate of machine \(m\) at period \(t\) under scenario \(s\).
    \item \( F^s_{m} \):= Failure time of machine \(m\) under scenario \(s\).
    \item \( D_{m,t} \):= Demand of product at machine \(m\) at period \(t\).

\end{itemize}

\subsection{Objective Function}
\begin{equation}
\label{eq:originalobj}
\begin{aligned}
\text{Minimize} \quad \frac{1}{|S|} \biggl[ 
&\underbrace{\sum_{s=1}^{S} \sum_{m=1}^{M} \sum_{t=1}^{F^s_m-1}  \left(C_{m,t} \cdot z_{m,t}\right)}_{\text{Preventive machine maintenance cost}} \\
&+ \underbrace{\sum_{s=1}^{S} \sum_{m=1}^{M} \sum_{t=F^s_m}^{T} \left( C_m^{f} \cdot z_{m,t} \right)}_{\text{Corrective machine maintenance cost}}\\
&+\underbrace{\sum_{s=1}^{S} \sum_{m=1}^{M} \left( \sum_{t=1}^{F^s_{m}-1} P^f \cdot z_{m,t} + \sum_{t=F^s_{m}}^{T} \left (t-F^s_{m}+1 \right ) \cdot P^f \cdot z_{m,t}\right)}_{\text{Penalty cost due to machine idle time}} \biggr]  \\
&+\underbrace{\sum_{m=1}^M \sum_{t=1}^{T}  \left(P^d \cdot d_{m,t}^{u}\right)}_{\text{Penalty cost due to unsatisfied demand}} \\
&+ \underbrace{\sum_{t=1}^T \delta_t \cdot \Delta}_{\text{Travel cost of the maintenance crew}}
\end{aligned}
\end{equation}

The objective function~\eqref{eq:originalobj} minimizes the total cost of operations and maintenance (O\&M) in leased manufacturing systems. It includes preventive and corrective maintenance costs, penalty costs associated with idle machine time and unmet demand, and travel costs of the maintenance crew. For preventive maintenance costs, we use the dynamic maintenance cost defined in \cite{karakaya2024sensor}. In this paper, the dynamic maintenance cost is defined as follows:

\begin{equation}
\label{DMC}
    C_{m,t} = \frac{C_m^{p} P(R_{t^o}^m > t) + C_m^{f} P(R_{t^o}^m \leq t)}{\int_0^\infty P(R_{t^o}^m > z) dz + t^o}
\end{equation}
where \(C_{m,t}\) is the sensor-driven estimate for conducting preventive maintenance at time \(t\) for machine \(m\), and \(R_{t^o}^m\) is the remaining life of machine \(m\) at observation time \(t^o\). The probabilities \(P(R_{t^o}^m > t\) and \(R_{t^o}^m \leq t)\) are driven by the predictions made using real-time sensor information. Here, \(C_m^{p}\) represents the planned (preventive) maintenance cost for the machine \(m \), and \(C_m^{f}\) denotes the failure replacement (corrective maintenance) cost for the machine \(m \), both of which are used to calculate the dynamic maintenance cost. For more details, one can refer to \cite{yildirim2016sensor}. In addition, we consider scenarios for machine failures. As a function of the failure time, we consider two cases. If the time is before failure, we use the dynamic maintenance cost function in \eqref{DMC}, as it balances the tradeoff between early and late maintenance decisions. However, for times after a failure, we impose a corrective cost (usually much larger than the preventive maintenance cost). 
These two cases correspond to the first and second lines of the objective function expression, respectively. 
We also have a penalty cost due to unsatisfied demands, idle machine costs reflect unproductive periods, and location change costs capture the expense of reallocating maintenance teams.

\subsection{Constraints}


\begin{equation}
\label{eq:lamda1}
\lambda^{s}_{m,t} \leq P^s_{m,t} \cdot (1- z_{m,t})  \quad \forall m \in M, \, t \leq F_{m}^s-1, \, s \in S
\end{equation}

\begin{equation}
\label{eq:lamda2}
\lambda^{s}_{m,t} \leq (P^s_{m,t} \cdot \sum_{l=1}^{t-1} z_{m,l}) \quad \forall m \in M, \, t \ge F_{m}^s, \, s \in S
\end{equation}


\begin{equation}
\label{eq:production_constraint1}
\gamma^{u,s}_{m,t} \geq D_{m,t}- \lambda^s_{m,t}  \quad \forall m \in M, \, t \in T, s \in \mathcal{S}
\end{equation}

\begin{equation}
\label{eq:production_constraint1_2}
d^{u}_{m,t} = \frac{1}{|S|}\sum_{s \in S} \gamma^{u,s}_{m,t}
\end{equation}

\begin{equation}
\label{eq:Maintenance Scheduling Constraints}
\sum_{m=1}^M z_{m,t} \leq J \quad \forall t \in T
\end{equation}

\begin{equation}
\label{eq:Max Maintenance}
\sum_{t=1}^T z_{m,t} = 1 \quad \forall m \in M
\end{equation}

\begin{equation}
\label{eq:Maintenance Location Consistency}
z_{m,t}  \leq \phi_{l,t} \quad \forall m \in M_n, \, l \in L, \, t \in T
\end{equation}

\begin{equation}
\label{eq:location changed}
 \phi_{l, t-1} \leq  \phi_{l, t} + \delta_t \quad \forall t \in \{2, \ldots, T\}, \, l \in L
\end{equation}

\begin{equation}
\label{eq:domain1}
\phi_{l,t} \in \{0, 1\} \quad \forall l \in L, \, t \in T
\end{equation}
\begin{equation}
\label{eq:domain2}
z_{m,t} \in \{0, 1\} \quad \forall m \in M, \, t \in T
\end{equation}
\begin{equation}
\label{eq:domain3}
\delta_{t} \in \{0, 1\} \quad \forall t \in T
\end{equation}

\begin{equation}
\label{eq:domain4}
d^{u}_{m,t} \geq 0 \quad \forall t \in T,  \, m \in M
\end{equation}

\begin{equation}
\label{eq:domain4_1}
\gamma^{u,s}_{m,t} \geq 0 \quad \forall t \in T,  m \in M, s \in {S}
\end{equation}

\begin{equation}
\label{eq:domain5}
\lambda_{m,t}^s \geq 0 \quad \forall t \in T,  \, s \in S, \, m \in M
\end{equation}

Constraints~\eqref{eq:lamda1} to \eqref{eq:production_constraint1_2} calculates the unsatisfied demand resulting from machine unavailability due to the failures or during maintenance period in the given scenario. We consider a first stage maintenance schedule that may either become a preventive or a corrective maintenance as a function of the failure time. If asset fails prior to maintenance, it becomes unavailable from the time of failure to the time of maintenance. Evidently, for any time periods after the failure time within a scenario, we check to see if a maintenance is scheduled prior to the time period. If not, it means that an asset has failed but was not yet fixed, thus we set the production to zero.
Constraint~\eqref{eq:Maintenance Scheduling Constraints} ensures that no more than $J$ maintenance activities are scheduled within a single period. Constraint~\eqref{eq:Max Maintenance} ensures conducting exactly one maintenance for each machine during planning horizon. Constraint~\eqref{eq:Maintenance Location Consistency} enforces that the maintenance team must be present at the location before maintenance begins. Constraint~\eqref{eq:location changed}  calculates whether the maintenance team changes location from one period to the next. Finally, constraints~\eqref{eq:domain1} to \eqref{eq:domain5} define the domains for the decision variables, ensuring that all binary and continuous variables adhere to their respective feasibility regions.

\section{MHA framework for operations and maintenance of leased manufacturing assets}
\label{MHA}
This section presents a multi-head attention transformation approach to tackle the operation and maintenance scheduling problem in leased manufacturing systems. In alignment with the MIP formulation, the proposed MHA transformation approach can also address the maintenance scheduling problem for leased manufacturing systems involving a set of \(M\) machines distributed across a set of \(N\) manufacturing sites, considered over a planning horizon of \(T\) time steps. Each machine is characterized by a combination of temporal and non-temporal features as we described in section~\ref{description} and \ref{MIP}.

\begin{figure}[h]
    \centering
    \includegraphics[width=\textwidth]{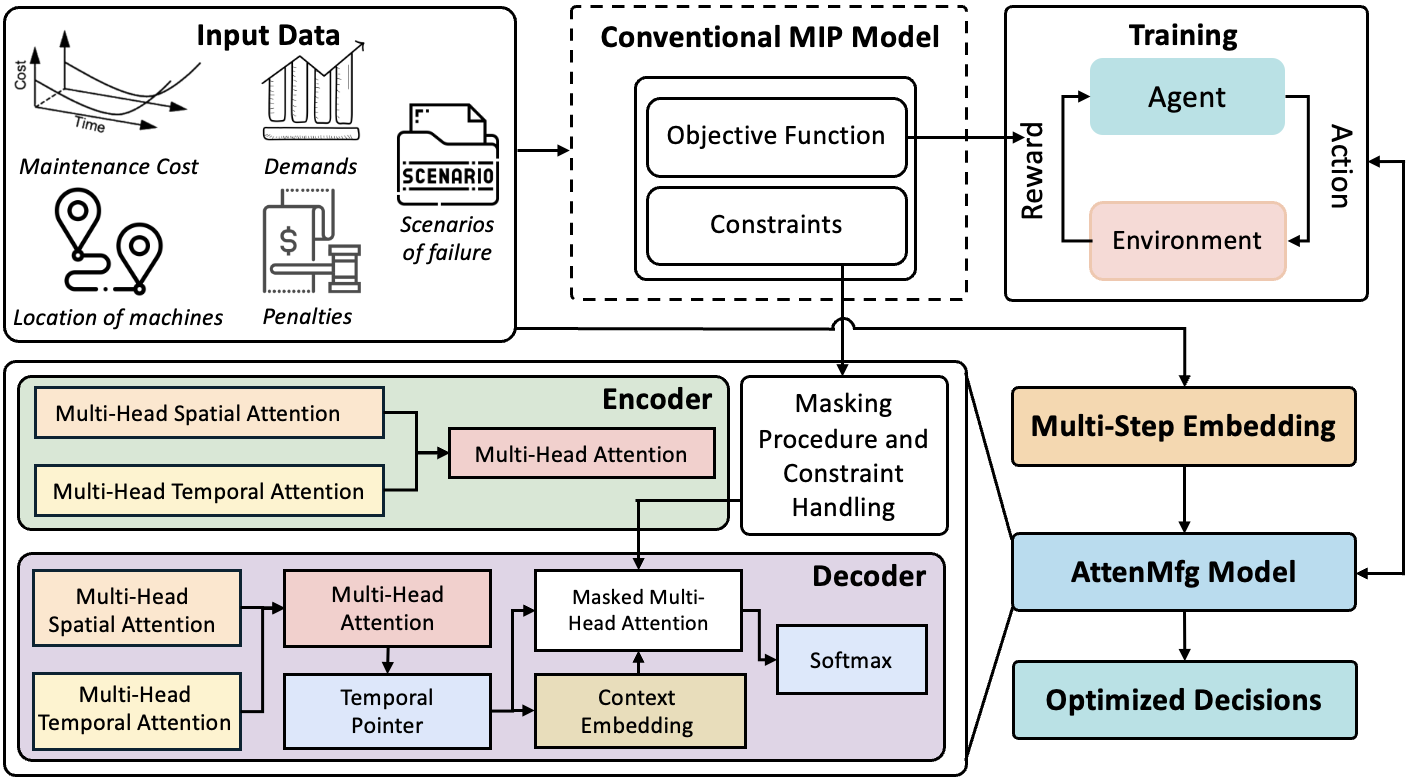}
    \caption{Framework for utilizing the proposed AttenMfg approach for operation and maintenance scheduling in leased manufacturing system. Input data, including dynamic maintenance costs, demands, penalties, and machine locations, is processed through multi-step embedding and encoded with multi-head spatial and temporal attention. The AttenMfg framework employs reinforcement learning for training and masking procedures to handle constraints, producing adaptive and efficient maintenance schedules.}
    \label{fig:how}
\end{figure}

To solve this problem using the MHA approach, we first need to embed the features discussed in Sections 3 and 4 into a representation compatible with the MHA methodology. Using this embedding, a reward function must then be defined to align with the objective function outlined in equation~\eqref{eq:originalobj}. Furthermore, a masking procedure is required to incorporate problem-specific constraints into the learning process, ensuring that the solutions generated are feasible. In the next subsection, we present a discussion of our embedding process that enables the MHA based learning.

\subsection{Multi-Step embedding for leased manufacturing maintenance problem}

To effectively solve the MIP problem in Section 4 using the MHA framework, it is essential to generate MHA-compatible embeddings for each of the maintenance and operation input features. 
The initial step in adopting an MHA based solution framework is to generate MHA-compatible embeddings for each of the maintenance and operation input features. 
For this purpose, we adopt the multi-step embedding approach proposed in \cite{kazemian2024attention}, as the concepts align closely with our requirements. 

\subsubsection{Maintenance cost embedding}
We incorporate preventive maintenance costs, failure costs, and failure scenarios as key input features for maintenance scheduling, which are used to generate embeddings compatible with MHA. Specifically, for machine \(m\) in period \(t\), the maintenance cost and associated penalties due to the unavailability are expressed as follows. When maintenance is performed before a failure occurs, the total cost includes the estimated maintenance cost and penalties caused by downtime during the maintenance period. In contrast, if a failure occurs before maintenance is performed, the total cost includes unexpected failure costs and penalty costs from the time of failure until maintenance is completed. Thus, the total maintenance scheduling cost, \(x_{m,t}\), for \(m \in M\) and \(t \in T\), can be formulated using equation~\eqref{eq:xit}:

\begin{equation}
\label{eq:xit}
\begin{aligned}
x_{m,t} = \frac{1}{|S|} \sum_{s \in S} \left[ I(t < F_m^s) \left( C_{m,t} + p^f \right) + I(t \geq F_m^s) \left( C_m^{pf} + \sum_{l=F_m^s}^{t} p^f \right) \right] 
\end{aligned}
\end{equation}

The maintenance cost embeddings derived from this formulation are critical for guiding the AttenMfg model in scheduling maintenance for a fleet of machines distributed in manufacturing sites under a defined set of scenarios.

\subsubsection{Throughput embedding}
In our MIP model, we incorporate a term for unsatisfied demand in the objective function as: \(
 \sum_{m=1}^M \sum_{t=1}^{T} \left(p^d \cdot d_{m,t}^{u}\right).
\) To model this, we use Algorithm 1, 
\begin{algorithm}
\caption{Finding Machine throughput ($\Lambda^{s,l}_{m,t}$)}
\label{alg:maintenance_parameters}
\begin{algorithmic}[1]
\FOR{all $m \in M$, $s \in S$}
        \STATE Find $F_m^s$
        \FOR{all $t \in T$}
            \IF{$t \leq F_m^s - 1$}
                \STATE $\Lambda^{s,l}_{m,t} = P_{m,t}^s (\forall l \neq t )$ 
                \STATE 
                $\Lambda^{s,t}_{m,t} = 0$ 
            \ELSIF{$t \geq F_m^s$}
                \STATE $\Lambda^{s,l}_{m,t} = 0$ ($\forall F_m^s \leq l \leq t $)
                \STATE $\Lambda^{s,l}_{m,t} = P_{m,t}^s$ ($\forall l > t $)
                \STATE $\Lambda^{s,l}_{m,t} = P_{m,t}^s$ ($\forall l < F_m^s $)
            \ENDIF
        \ENDFOR
\ENDFOR
\end{algorithmic}
\end{algorithm}
which outlines the process of determining the production limits for each machine (\(m\)) under different failure time scenarios (\(s\)) and maintenance schedules. First, we determine the failure time for each machine in a given scenario (\(F_m^s\)). Based on this failure time, we define the production limits of machine \(m\) under scenario \(s\) at time period \(l\) when we conduct maintenance at time period \(t\), (\(\Lambda^{s,l}_{m,t}\)),  are set either to pre-failure limits and post-failure limits as algorithm~\ref{alg:maintenance_parameters}. In equation~\eqref{eq:throughput}, based on the average production level for each machine at each time, we determine the unsatisfied demand penalties.

\begin{equation}
\label{eq:throughput}
y_{m,t} =  \left( p^d \cdot  \frac{1}{|S|}\sum_{s \in S} \sum_{l \in T} \left\{ \max\left(D_{m,l} -  \Lambda^{s,l}_{m,t},0\right) \right\} \right),
\end{equation}

By embedding \( y_{m,t} \), we represent the penalty cost as part of the throughput embedding in the input feature.

\subsubsection{Idle period embedding}
In the multi-head attention-based model, at each timestep, a machine must be selected for maintenance, ensuring no idle periods; however, constraint~\eqref{eq:Maintenance Scheduling Constraints} in our MIP allows for the possibility of idle periods, where no maintenance is performed. We model idle periods by modeling a visit to an idle machine, which corresponds to a no maintenance period. More specifically, we define a set of idle machines (\( M'\)) with constraints such that $x_{m,t} = 0, y_{m,t} = 0 \quad \forall m \in M', \forall t$.
These constraints ensure that machines in the set \( M' \) incur no maintenance cost while stationed at the depot or base. These equations imply that the crew's location coincides with the depot whenever maintenance for idle machines is scheduled. As a result, the set of idle machines helps MHA model relax the requirement of exactly \( J \) decisions per period.

\subsubsection{Crew logistics embedding}
\label{Crew logistics embedding}
As mentioned earlier in the context of idle period embedding, the MHA model selects only one machine for maintenance at each timestep. However, with constraint~\eqref{eq:Maintenance Scheduling Constraints}, our MIP allows up to \(J\) opportunistic maintenance operations per period. To accommodate this, we duplicate the input features, \(x_{m,t}\) and \(y_{m,t}\) representing the maintenance scheduling cost for machine \(i\) and penalties due to unsatisfied demands at period \(t\), based on the number of allowable maintenance operations per period. As a result, we obtain transformations that ensure that $x \in R(M \times T) \longrightarrow{\chi \in R(M \times T \times J)}$ and $y\in R(M \times T)  \longrightarrow{Y\in R(M \times T \times J)}$.


\subsubsection{Dimension alignment embedding}
This step ensures that the dimensional consistency of manufacturing site location data aligns with the features used for maintenance decision costs and throughputs. Initially, the manufacturing site location data \(N\) for each machine is represented in \( R(M+M') \), while \( \chi \) and \(Y\) exists in \( R(M+M') \times T \times J \). To harmonize these dimensions, the location vector \( M+M' \) is repeated \( T \times J \) times, results in \(N'\) aligning it with the dimensionality of \( \chi \text{ and } Y\). This alignment integrates location information with other input features, maintaining uniformity across the dataset. We therefore enable a transform wherein $N \in R(M+M') \longrightarrow{N' \in R(M+M') \times T \times J}$.


\subsubsection{Input Embedding} The \emph{Input Embedding} layer maps the inputs (\( \chi, Y, N' \)) into a unified \( D \)-dimensional space. This transformation uncovers patterns and dependencies, enabling more effective maintenance scheduling and improving the model's overall performance.

\begin{equation}
\label{eq:Input Embedding}
E_{\chi, Y, N'} = \text{MHA Embedding}(\chi, Y, N')
\end{equation}

\subsection{Cost Function}
Our attention-based model uses a cost function that corresponds directly to the objective function outlined in the MIP model \eqref{eq:originalobj}. This correspondence allows our deep learning model to optimize the same cost parameters, effectively bridging the gap between traditional optimization frameworks and advanced machine learning approaches in maintenance planning. The cost of performing maintenance on machine \(v\) at time \(t\), following maintenance on machine \(u\), is defined by the function 
\(f_c : (\chi_{u, t-1} + Y_{u,t-1}) \times (\chi_{v, t}+ Y_{v,t}) \to \mathbb{R}\). 
The objective is to identify a sequence \(\pi\) of machines with length  \(\in (0, T \times J]\) that meets constraints~\eqref{eq:Maintenance Scheduling Constraints} and \eqref{eq:Max Maintenance}, as verified by \(C(\pi)\). Here, \(T \times J\) accounts for extended time periods resulting from repeated maintenance actions within each period (\ref{Crew logistics embedding}). The optimization problem aims to minimize:
\begin{equation}  
P_{\text{obj}}[\pi \mid (M+M')] = \sum_{t=1}^{T \times J} f_c(\chi(t-1, u)+Y(t-1, u), \chi(t, v) +Y(t, v))
\end{equation}

where \(u, v \in I\) represent machines chosen from the set \(M+M'\) at time \(t\). The constraint-checking function \(C\) ensures that sequence \(\pi\) satisfies the required conditions.

The maintenance cost function \(f_c(\chi(t-1, u)+Y(t-1, u), \chi(t, v) +Y(t, v))\) is defined as:
\begin{equation}  
f_c(\chi(t-1, u)+Y(t-1, u), \chi(t, v)+Y(t, v)) = \chi(t, v) +Y(t, v) + \mathbb{I}(\delta_t = 1) \cdot \Delta
\end{equation}
where \(\mathbb{I}(\delta_t = 1)\) is an indicator function that activates when \(\delta_t = 1\), applying an additional cost \(\Delta\) if the location of maintenance team change in \(t\).

\section{Attention Model}
\label{attention}
The objective of the leased manufacturing machine maintenance scheduling problem is to determine a policy \(\Pi\) that generates a sequence \(\pi\), minimizing the objective \(P_{\text{obj}}\) while satisfying the constraint \(C(\pi)\) for a given problem instance \(s\). The proposed framework, AttenMfg, leverages a MHA mechanism and is illustrated in figure~\ref{fig:how}. It consists of an encoder-decoder structure that defines a stochastic policy \(p(\Pi | s)\), parameterized by \(\theta\), as shown in equation~\ref{eq:pi_cost}:

\begin{equation}
\label{eq:pi_cost}
p_{\theta}(\pi | s) = \prod_{t=1}^{T \times J} p_{\theta}(\pi_t | s, \pi_{1:t-1}).
\end{equation}

\subsection{Encoder}

The encoder employs a Graph Temporal Attention (GTA) network inspired by \cite{gunarathna2022solving,kazemian2024attention}, built upon the transformer architecture from \cite{vaswani2017attention}. It processes inputs through parallel spatial and temporal attention layers, combining their outputs (\(H_S^{(1)}, H_T^{(1)}\)) to produce a comprehensive sequence embedding. Additional layers refine these embeddings, where the hidden state of machine \(m\) at time \(t\) in layer \(l\) is denoted as \(h_{m,t}^{(l)}\).

\subsubsection{Spatial Attention Layer}

The spatial attention layer models inter-machine dependencies at a given time step. For a given time \(t\), the sequence of machine features is represented as \(H_{S,t}^{(l)} = \{h_{1,t}^{(l)}, h_{2,t}^{(l)}, \ldots, h_{I,t}^{(l)}\} \in \mathbb{R}^{(M+M') \times D^h}\). Using self-attention \cite{vaswani2017attention}, the query, key, and value projections are computed:

\begin{equation}
q_{t}^{S(l)} = w_{q}^{S(l)} H_{S,t}^{(l)}, \quad k_{t}^{S(l)} = w_{k}^{S(l)} H_{S,t}^{(l)}, \quad v_{t}^{S(l)} = w_{v}^{S(l)} H_{S,t}^{(l)}.
\end{equation}

The attention weights are calculated as:

\begin{equation}
\alpha_t = \text{softmax}\left( \frac{q_{t}^{S(l)} (k_{t}^{S(l)})^T}{\sqrt{D^h}} \right),
\end{equation}

and the spatial attention output is:

\begin{equation}
H_{S,t}^{(l+1)} = \alpha_t v_{t}^{S(l)}.
\end{equation}

\subsubsection{Temporal Attention Layer}

The temporal attention layer captures dependencies between machine features across different time steps. For input \(H^{(l)}_{T}\), the query, key, and value projections are:

\begin{equation}
q_{m}^{T(l)} = w_q^{T(l)} H_{T,m}^{(l)}, \quad k_{m}^{T(l)} = w_k^{T(l)} H_{T,m}^{(l)}, \quad v_{m}^{T(l)} = w_v^{T(l)} H_{T,m}^{(l)}.
\end{equation}

The attention weights and output are computed as:

\begin{equation}
\beta_m = \text{softmax}\left( \frac{q_m^{T(l)} (k_m^{T(l)})^T}{\sqrt{D^h}} \right),
\end{equation}

\begin{equation}
H_{T,m}^{(l+1)} = \beta_t v_{m}^{S(l)}.
\end{equation}

\subsubsection{Integration Layer}

Spatial and temporal outputs (\(H_S^{(l+1)}, H_T^{(l+1)}\)) are concatenated, linearly transformed, and passed through a sigmoid activation:

\begin{equation}
H^{(l+1)} = \sigma(wI \cdot (H_S^{(l+1)} || H_{T}^{(l+1)})).
\end{equation}

\subsection{Decoder}

The decoder uses the encoder output and iteratively constructs the solution sequence. At each step, it generates embeddings for unmaintained machines, combines them with contextual information, and selects the next machine using a masked multi-head attention layer.

\subsubsection{Temporal Pointer}

The temporal pointer dynamically computes attention weights at each decoding step. Slicing the encoder embedding at time \(t\), it focuses on the most relevant information. The operations are defined as in \cite{gunarathna2022solving}.

\subsubsection{Context Embedding}

The context embedding, \(H_C\), incorporates problem-specific information, including the last maintained machine and aggregated features of all machines at time \(t\):

\begin{equation}
H_{C,t} = \{ h^{(L)}_{\pi_t} \| H_{G,t} \}; \quad H_{G,t} = \sum_{m=0}^{M+M'} h^{(L)}_{m,t}.
\end{equation}

\subsubsection{Masked Multi-head Attention and Log-Probability Layer}

The masked MHA layer integrates the context-specific embedding \(H_{C,t}\) with the current problem representation \(H_{D,t}\) to compute probabilities for selecting the next machine. This is achieved by projecting the query, key, and value weights as follows:

\begin{equation}
qC_t = w_{C_q} H_{C,t}, \quad kC_t = w_{C_k} H_{D,t}, \quad vC_t = w_{C_v} H_{D,t}.
\end{equation}

The attention weights are computed using the scaled dot-product mechanism:

\begin{equation}
\alpha_t = \text{softmax}\left( \frac{qC_t kC_t^T}{\sqrt{D^h}} \right),
\end{equation}

and the resulting embedding \(H_{D,t}^{(F)}\) is:

\begin{equation}
H_{D,t}^{(F)} = \alpha_t vC_t.
\end{equation}

To compute the log-probabilities for each machine, a tanh activation function and a weight vector \(w_P\) are applied:

\begin{equation}
\gamma_t = \tanh\left( H_{D,t}^{(F)} \cdot (w_P \cdot H_{D,t})^T \right).
\end{equation}

Finally, the log probabilities for selecting the next machine are calculated through a softmax layer:

\begin{equation}
P_t = \text{softmax}(\gamma_t).
\end{equation}

\paragraph{Masking Procedure:} 
To ensure compliance with constraints, a masking mechanism is applied to \(P_t\). Machines that have already been maintained or whose selection would violate constraints are assigned a large negative value in their probability score. This is achieved by updating the mask \(M\) at each time step. The masked probabilities are computed as:

\begin{equation}
P_t = \text{softmax}(\gamma_t + M),
\end{equation}

where \(M\) is initialized as a zero vector and updated based on maintenance constraints. If a machine selection violates a constraint, \(M[i] = -\infty\) is applied to exclude it from the next selection.

This mechanism ensures that the selection of machines adheres to all constraints while dynamically adjusting probabilities based on the current state of the solution.

\paragraph{Constraint handling:} Unlike traditional MIP models, where adding such constraints can significantly increase complexity—making the problem harder to solve or even unsolvable—our approach maintains efficiency and scalability. Our attention-based maintenance scheduling model seamlessly incorporates various operational constraints, ensuring both feasibility and efficiency. These constraints include limits on preventive maintenance actions per period, total allowable maintenance tasks, crew availability, and logistical restrictions. This is achieved through a dynamic masking process in the decoder, a feature that enables the model to exclude actions violating constraints by setting their attention weights to zero.

\subsection{Reinforcement Learning-driven AttenMfg}

Our model leverages a reinforcement learning approach for training, where the probability distribution \(p_{\theta}(\pi|s)\) defines a solution \(\pi\) for a given instance \(s\). The loss function, \(L(\theta|s)\), represents the expected cost, combining maintenance expenses, unsatisfied demand penalties and team relocation costs. Using gradient descent, the reinforcement learning gradient estimator is employed to optimize the model parameters, as shown below:

\begin{equation}
\nabla_{\theta}L(\theta|s) = \mathbb{E}_{p_{\theta}(\pi|s)} \left[ \left( L(\pi) - b(s) \right) \nabla_{\theta} \log p_{\theta}(\pi|s) \right]
\end{equation}

To stabilize training and reduce gradient variance, a baseline function \(b(s)\) is introduced. This reference point aids in accelerating learning by providing a cost comparison. We use a rollout-based baseline, where the average cost is calculated by simulating the policy multiple times. For each state \(s\), \(b(s)\) approximates the expected cost of following the current policy by averaging costs across these rollouts. This reinforcement learning-driven approach ensures efficient training and robust policy optimization.

\subsubsection{Training dataset types}

To evaluate the performance AttenMfg for solving operations and maintenance in leased manufacturing systems, we developed a set of diverse training model types. These model types incorporate critical factors such as dynamic maintenance costs, failure costs, and failure time scenariosdrawing on methodologies outlined in \citep{gebraeel2005residual,yildirim2016sensor1}. These parameters are derived from real-world data collected from a rotary machine, ensuring that the training dataset reflects a realistic operational environment to enhance its practical relevance and applicability.

Additionally, the demands were generated using a random process. For each instance, the demands \(D_m,t\) were sampled from a normal distribution, \(D_m,t \sim \mathcal{N}(\mu, 0.1 \times \text{max production limit})\), where the mean \(\mu\) was set equal to the maximum production limit.

\hspace{-6mm}\textit{Machine Location Assignment:} The allocation of machines to manufacturing sites in our datasets was addressed using two distinct types:

\begin{itemize}
    \item \textit{Type 1-Fixed Number of Sites:}  
    In this type, a fixed number of manufacturing sites was predefined. Each machine was then randomly assigned to one of these sites. This approach maintained a consistent number of sites across all dataset instances, allowing for controlled evaluation of the model's performance.
    
    \item \textit{Type 2-Variable Number of Sites:}  
    In this type, the number of manufacturing sites was dynamically determined for each dataset instance. Specifically, the number of sites was randomly selected between 1 and 10. Once the number of sites was established, machines were randomly assigned to these sites, ensuring variability in their distribution. This method introduced diversity into the datasets, simulating a wide range of operational scenarios and enhancing the model's adaptability to different real-world conditions.
\end{itemize}

\section{Experiments}
\label{experiment}

In this study, we propose an attention-based transformation model to address the operations and maintenance (O\&M) scheduling problem in leased manufacturing systems. The experimental analysis is structured into three parts to evaluate the model's performance comprehensively.

In the first case, we aim to showcase the performance of our model and architecture in terms of computational time and solution quality. To achieve this, we compare our model against a benchmark method, which involves solving the Mixed-Integer Programming (MIP) formulation of the problem (introduced in Section 4) using a conventional commercial solver, such as Gurobi. This comparison highlights the efficiency and effectiveness of our attention-based approach relative to traditional optimization techniques. The second case focuses on assessing the generalizability of our model when tested on datasets with a fixed number of sites. Specifically, we examine whether a model trained on one dataset instance with a predefined number of manufacturing sites can effectively solve problems from other dataset instances. The third case investigates the generalizability of our model when the number of manufacturing sites varies across a dataset instances. In this analysis, we evaluate whether a model trained on datasets with a specific distribution of variable site numbers can effectively handle problems from other datasets where the site numbers follow a different distribution or configuration. By incorporating this variability, we aim to assess the model's robustness and adaptability to more diverse and dynamic real-world scenarios. This analysis highlights the impact of introducing randomness in the number of manufacturing sites on the model's generalizability.

To facilitate the experiments, we define distinct problem settings, summarized in Table~\ref{tab:problem_settings_rotated}. The table provides an overview of the problem configurations used for different cases. These settings are denoted by our model naming convention \(LxPxMx\), where \(Lx\) represents the number of manufacturing sites, \(Px\) the number of planning periods, and \(Mx\) the number of machines. For Type~2 models, the number of sites is denoted by \(R\), indicating that it is randomly determined for each instance. For example, \(LRP15M40\) refers to a model with a random number of sites, 15 planning periods, and 40 machines. In all configurations, the total number of maintenance teams available per period is fixed at 3. In the following, we use this naming convention not only to refer to the models but also to the datasets corresponding to each configuration. This dual use ensures clarity when discussing experimental setups and results.

\begin{table}[h!]
\centering
\footnotesize
\caption{Model Settings}
\label{tab:problem_settings_rotated}
\renewcommand{\arraystretch}{1.3} 
\setlength{\tabcolsep}{8pt} 
\begin{tabular}{@{}lcccccc@{}}
\toprule
\textbf{Model Type} & \textbf{Model Name} & \textbf{\# of Sites} & \textbf{\# of Periods} & \textbf{\# of Machines} \\ \midrule
\multirow{6}{*}{Type 1} 
& \textbf{L5P10M25}   & 5                    & 10                     & 25                      \\
& \textbf{L5P15M40}   & 5                    & 15                     & 40                      \\
& \textbf{L5P20M50}   & 5                    & 20                     & 50                      \\
& \textbf{L10P10M25}  & 10                   & 10                     & 25                      \\
& \textbf{L10P15M40}  & 10                   & 15                     & 40                      \\
& \textbf{L10P20M50}  & 10                   & 20                     & 50                      \\ \midrule
\multirow{2}{*}{Type 2} 
& \textbf{LRP15M40}   & R*                   & 15                     & 40                      \\
& \textbf{LRP20M50}   & R*                   & 20                     & 50                      \\ \bottomrule
\end{tabular}
\bigskip
\captionsetup{font=small}
\caption*{\small R*: Random number between 1 and 10.}
\end{table}

The proposed model uses $N=3$ encoder layers with a hidden dimension of $D_h=128$, balancing solution quality and computational complexity. Training is performed over 100 epochs with 12,800 problem instances, using a batch size of 16 and a learning rate of $1 \times 10^{-4}$ to ensure stability. The Adam optimizer is applied for effective convergence. Experiments were conducted on a supercomputer with an Intel processor and two NVIDIA Quadro RTX 6000 GPUs, providing the computational resources needed for efficient training and evaluation.

\subsection{Case Study 1: Comparative Performance of Proposed and Benchmark Method}

\begin{figure}[h]
     \centering
     \includegraphics[width=0.65\linewidth]{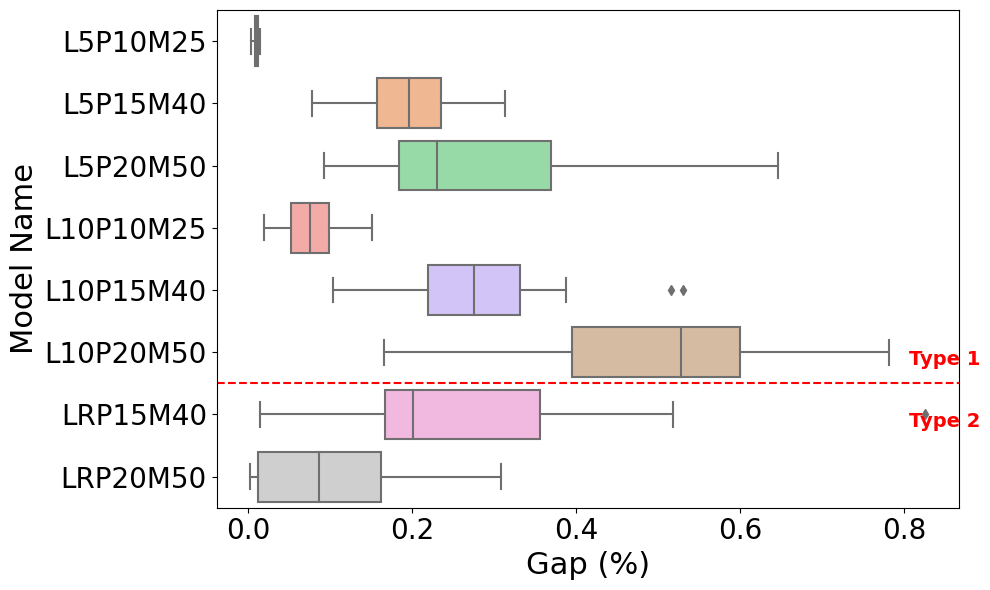} 
     \caption{Summary of Solution Gaps for 20 Problems in Each Model.}
     \label{fig:results_summary}
 \end{figure}

In this case, we solved 20 problem instances from each dataset configuration. Each problem instance was first solved using the Mixed Integer Programming (MIP) approach with Gurobi to obtain the exact optimal solution. Subsequently, the same 20 problem instances were solved using our proposed attention-based method, AttenMfg. The performance of AttenMfg was then evaluated by calculating the optimality gap (\%) between the exact solutions obtained from Gurobi and the solutions generated by AttenMfg. Figure~\ref{fig:results_summary} provides a summary of the results, presenting both the solution gaps and their distributions across the problem instances. This comparison highlights the effectiveness and accuracy of our method relative to the exact benchmark solutions.

Solving these problems using a Mixed-Integer Programming (MIP) approach with a commercial solver proved highly time-consuming (Table~\ref{table:solving_time}). Optimal solutions were found for cases \(L5P10M25, L10P10M25\) and \(L5P15M40\) within 2 hours. However, for other cases, including larger and more complex configurations, the solver could not find optimal solutions within the same time limit. In contrast, our method solved all problem instances across all settings in less than a second, demonstrating significant computational efficiency. In terms of solution quality, our approach achieved a mean gap of less than 0.2\% for the smaller configurations, such as \(L5P10M25\) and \(L10P10M25\). For larger instances, such as \(L10P20M50\) (Type 1) and \(LRP20M50\) (Type 2), the mean gap remained within 0.4–0.5\%, demonstrating the effectiveness of our model even as the problem complexity increased. Notably, the solution gaps for Type~2 instances (with a variable number of sites) still maintained a mean gap below 0.5\%, highlighting the robustness of our approach across diverse configurations. The box plot highlights the distribution of gaps, showing tight distributions in simpler settings and slightly larger variability in more complex ones. These results demonstrate the effectiveness of our model in both solution time and quality.


\begin{table}[ht]
    \centering
    \caption{Solution times for MIP model across different cases. Corresponding solution times for AttenMfg were excluded, as it solves all cases under one second.}
    \begin{tabular}{@{}lccccc@{}}
        \toprule
        Model Name  & Mean    & Median  & Q1      & Q3      & \% Solved \\ \midrule
        L5P10M25    & 486.25  & 495.00  & 401.25  & 552.50  & 100\%     \\
        L5P15M40    & 989.80  & 815.00  & 609.00  & 1250.00 & 100\%     \\
        L5P20M50    & 7200.00 & 7200.00 & 7200.00 & 7200.00 & 0\%     \\
        L10P10M25   & 855.10  & 855.50  & 658.50  & 980.25  & 100\%     \\
        L10P15M40   & 7200.00 & 7200.00 & 7200.00 & 7200.00 & 0\%     \\
        L10P20M50   & 7200.00 & 7200.00 & 7200.00 & 7200.00 & 0\%     \\
        LRP15M40    & 7200.00 & 7200.00 & 7200.00 & 7200.00 & 0\%     \\
        LRP20M50    & 7200.00 & 7200.00 & 7200.00 & 7200.00 & 0\%     \\ \bottomrule
    \end{tabular}
    \label{table:solving_time}
\end{table}

\subsection{Case Study 2: Generalizability Over Fixed Number of Sites (Type 1)}

\begin{figure}[htbp]
    \centering
    \begin{subfigure}{0.30\textwidth}
        \centering
        \includegraphics[width=\linewidth]{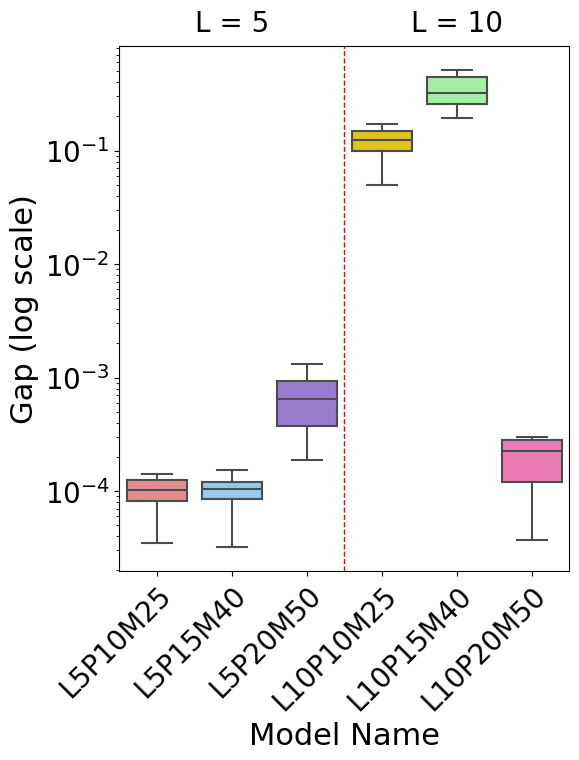}
        \caption{Evaluating on L5P10M25}
    \end{subfigure}
    \hfill
    \begin{subfigure}{0.30\textwidth}
        \centering
        \includegraphics[width=\linewidth]{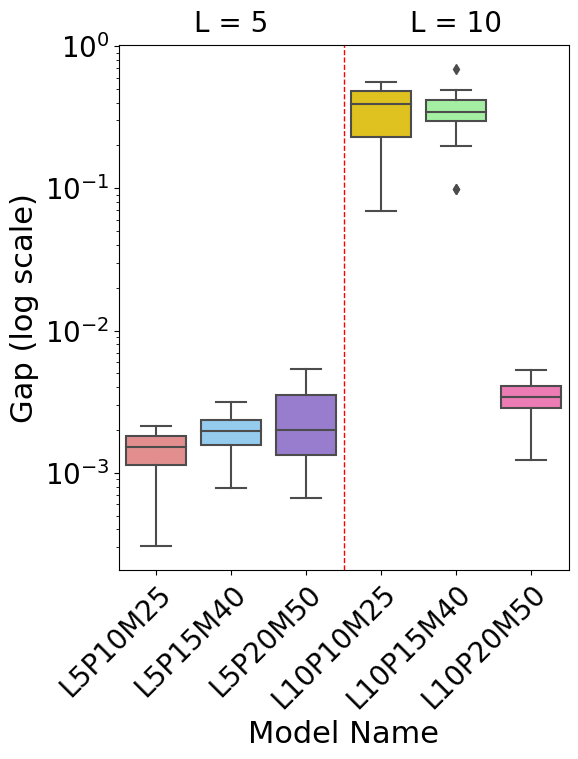}
        \caption{Evaluating on L5P15M40}
    \end{subfigure}
    \hfill
    \begin{subfigure}{0.30\textwidth}
        \centering
        \includegraphics[width=\linewidth]{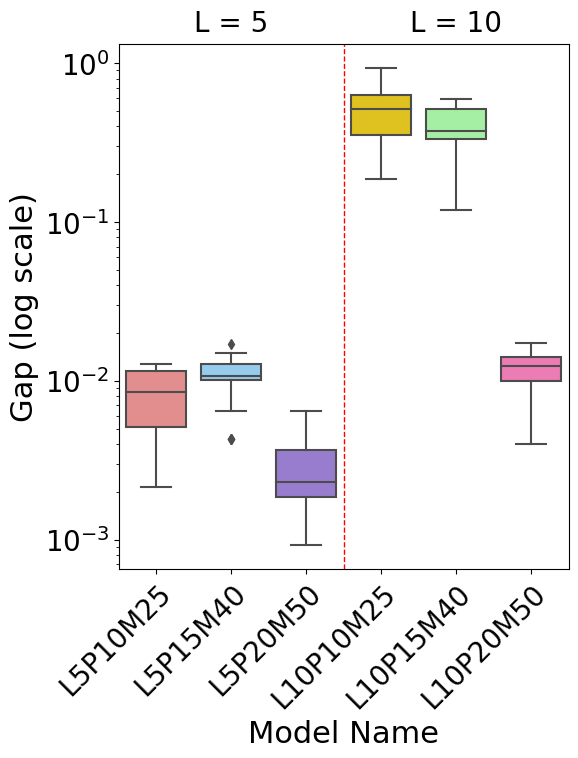}
        \caption{Evaluating on L5P20M50}
    \end{subfigure}
    \vspace{0.2cm}
    \begin{subfigure}{0.30\textwidth}
        \centering
        \includegraphics[width=\linewidth]{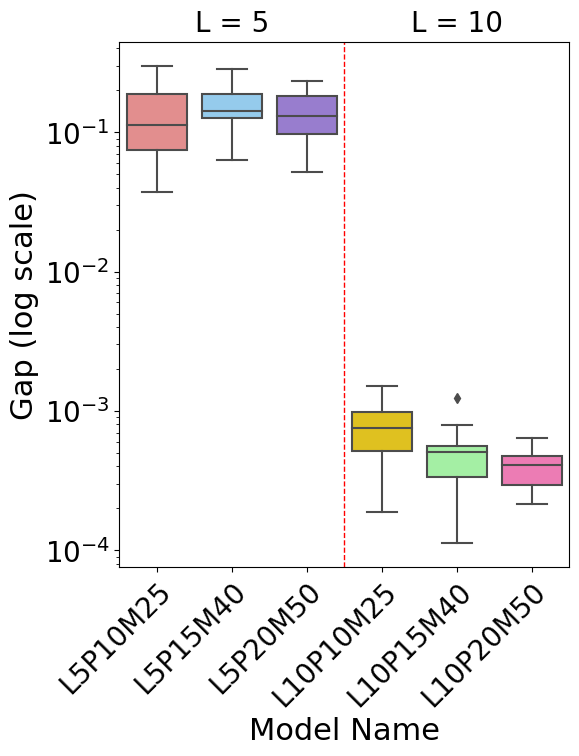}
        \caption{Evaluating on L10P10M25}
    \end{subfigure}
    \hfill
    \begin{subfigure}{0.30\textwidth}
        \centering
        \includegraphics[width=\linewidth]{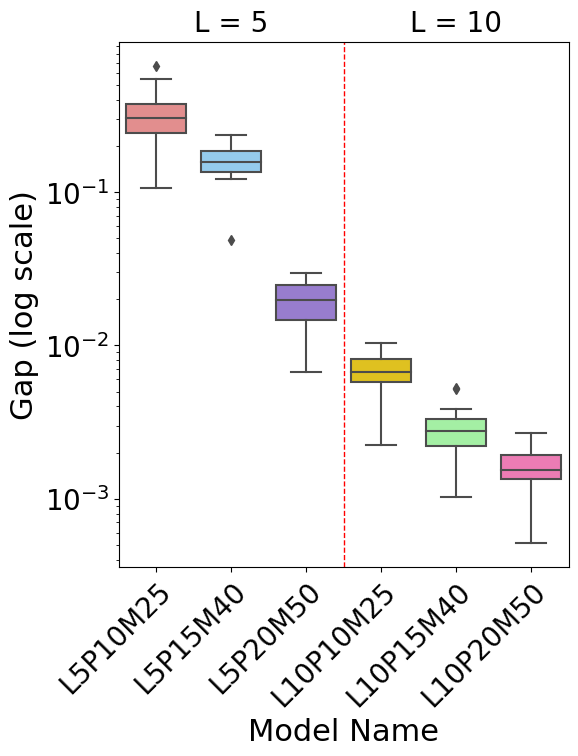}
        \caption{Evaluating on L10P15M40}
    \end{subfigure}
    \hfill
    \begin{subfigure}{0.30\textwidth}
        \centering
        \includegraphics[width=\linewidth]{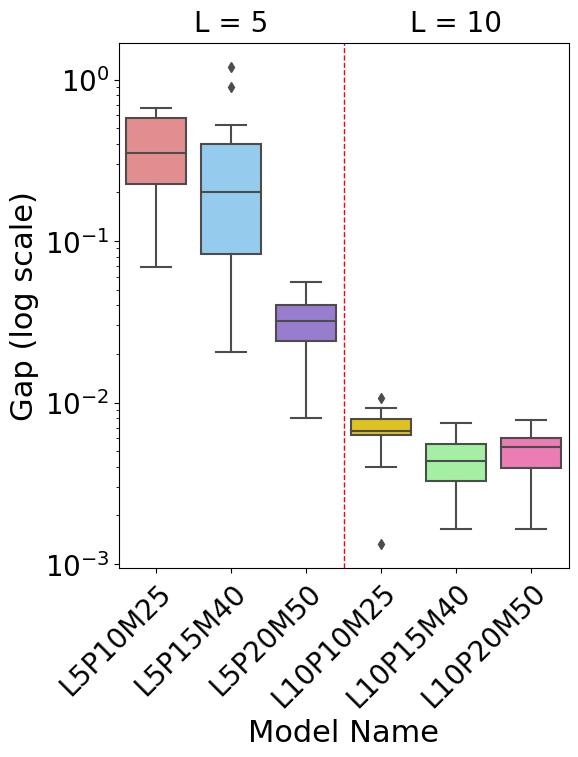}
        \caption{Evaluating on L10P20M50}
    \end{subfigure}

    \caption{Generalizability analysis: solution gaps (log scale) when models trained on different datasets (x-axis) are evaluated on fixed target datasets (a–f). The results highlight that models perform best when the number of sites in the training and testing datasets match, while transitioning to larger configurations with different site numbers results in a significant performance drop, with a maximum observed gap of 53.09\%. Additionally, models trained on larger and more complex datasets demonstrate strong backward generalizability, improving performance when evaluated on smaller datasets.}

    \label{fig:combined}
\end{figure}

In this case, we analyze the generalizability of our AttenMfg model by evaluating its performance when trained on one dataset configuration and tested on other configurations. Specifically, for each subfigure in figure~\ref{fig:combined}, the model is trained on the dataset corresponding to the model name shown on the x-axis and evaluated on a specific target dataset, as indicated in the subfigure titles. The solution quality is measured using the optimality gap (log scale) between the solutions generated by our model and the exact solutions.

From Figure~\ref{fig:combined}, several key trends emerge:
\begin{itemize}
    \item \textit{Accuracy is higher while generalizing for cases with matching site numbers:} As observed in Figure~\ref{fig:combined}, evaluating the model on a dataset with the same number of manufacturing sites as in the training dataset consistently yields better results. This indicates that the model performs more accurately when the training and testing configurations align in terms of the number of sites.
    
    \item \textit{It is more challenging to generalize to larger cases with different site numbers:} When the model is trained using small cases, evaluating the model on larger cases with different number of sites results in a significant performance drop, with a maximum observed gap of 53.09\%. This highlights the difficulty of generalizing effectively when scaling up both the problem size and the number of sites, underscoring the model's sensitivity to such changes.

    \item \textit{High accuracy is attained during backward generalizability:} The model demonstrates strong backward generalizability, meaning that training on larger and more complex datasets improves performance when evaluated on smaller datasets. This trend highlights the benefit of training on datasets with greater complexity to enhance the model's robustness across varying problem sizes.
\end{itemize}

These findings offer significant managerial implications for leased manufacturing systems: \textit{Lessors do not have to train their models for each particular case. Instead, a good practical approach is to train a single model using a problem setting that reflects the most complex operational scenario, and use it while solving all the cases.}

In conclusion, the ability to generalize across different problem settings highlights the adaptibility of our AttenMfg model. By strategically selecting training settings, lessors can effectively solve current and future problems, making this approach a valuable tool for long-term operational planning.

\subsection{Case Study 3: Generalizability with Variable Number of Sites (Type 2)}
In this case study, we assess the generalizability of our AttenMfg model when trained on datasets with a random number of sites (Type 2) and evaluated on all other configurations, including both fixed (Type 1) and random site settings. The objective is to investigate whether training on datasets with greater variability in the number of sites improves the model's robustness and adaptability to diverse problem settings. 

\begin{figure}[h]
     \centering
     \includegraphics[width=0.7\linewidth]{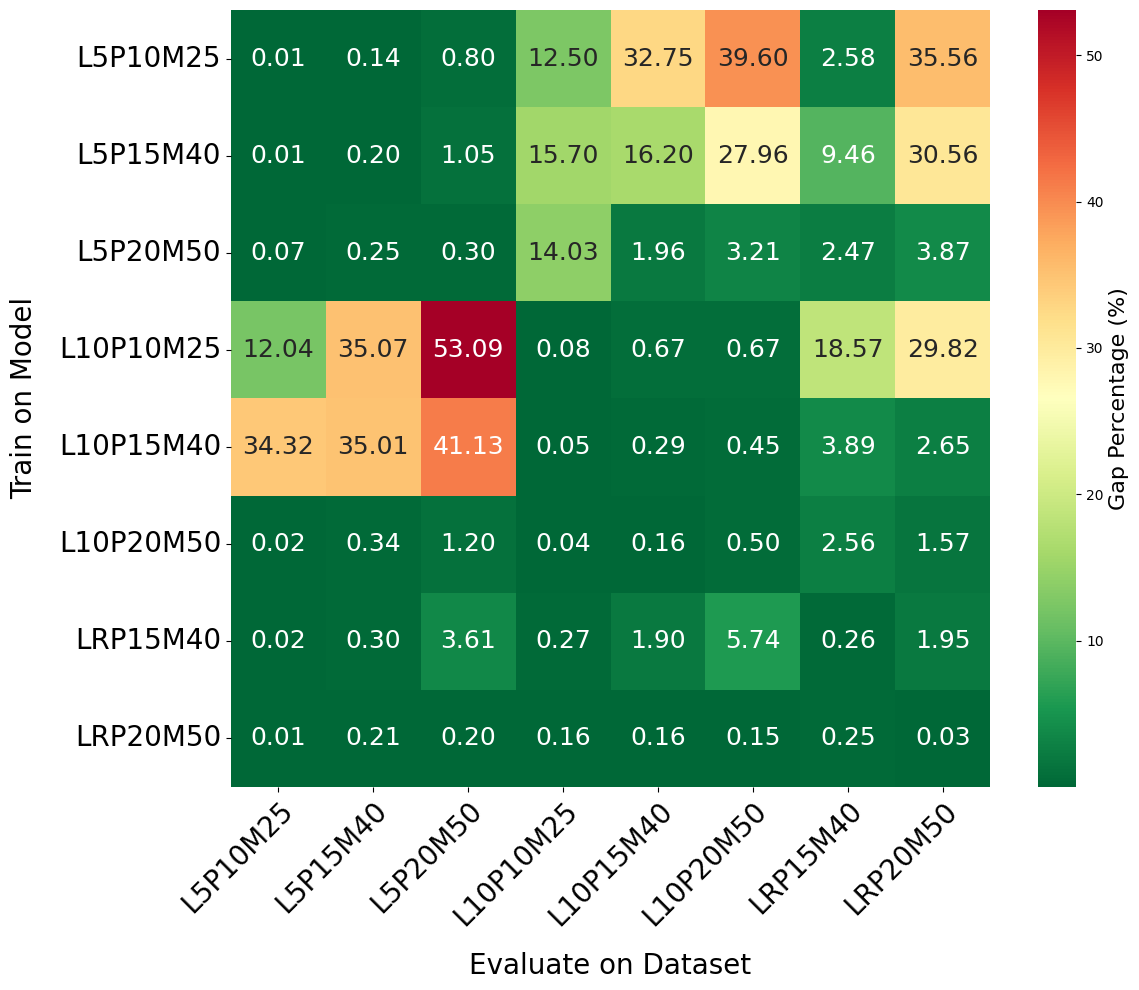} 
     \caption{Heatmap of gap percentages: models trained on random-site datasets evaluated across all configurations. The results highlight the strong generalization performance of random-site models (\(LRP15M40\) and \(LRP20M50\)), achieving low solution gaps across diverse configurations. Notably, the larger random model (\(LRP20M50\)) performs comparably to testing on the same model and constructing new datasets for specific configurations.}

     \label{fig:results_summaryRR}
 \end{figure}

The results are summarized in Figure~\ref{fig:results_summaryRR}, which presents the optimality gap percentage for each evaluation scenario. The y-axis represents the datasets on which the model was trained, and the x-axis corresponds to the datasets on which the model was evaluated. Key findings from the results are as follows:

\begin{itemize}
    \item \textit{Training on random cases yield strong generalization across all problems:} Training on cases with a random number of sites (e.g., \(LRP15M40\) and \(LRP20M50\)) yields strong generalization performance across all problem settings. The models trained on these random-site datasets achieve low optimality gaps when evaluated on fixed-site configurations, highlighting the benefit of incorporating variability in the training datasets to enhance the model's robustness and adaptability.
    \item \textit{Larger models trained on random cases improves adaptibility \& solution quality:} The larger random model (\(LRP20M50\)) demonstrates performance comparable to the models that are trained and tested in the same particular case, which are very computationally expensive as it would require training a new model for each case. Instead, a single large randomized model can generalize with a similar solution quality across a diverse set of problem settings,  minimizing the need for customized datasets for each model.

\end{itemize}

\subsection{Managerial Insights}
The extensive case studies yield several actionable managerial insights for lessors to enhance operational efficiency and maintenance decision-making in leased manufacturing systems. These insights can be summarized as follows:

\textit{Sensor-driven proactive maintenance management is key to the success of leased manufacturing systems:}  
    Integrating sensor-driven insights within an attention-based framework enables lessors to dynamically optimize maintenance schedules as a function of the sensor-driven insights on the condition of the machines. By minimizing machine downtime, unsatisfied demand penalties, and relocation costs, this approach not only enhances operational efficiency but also supports sustainability goals by reducing unnecessary maintenance actions and eliminating failure events.

    \textit{The proposed model inherently addresses the lack of adaptability and high computation costs in existing models:} 
    A significant barrier to adopting advanced maintenance scheduling technologies is their lack of adaptability to diverse configurations (see Figure \ref{fig:why}). Existing methods require constant engineer supervision to adapt to different problem settings, which proves costly especially for enterprises with limited access to resources. This study addresses this challenge by developing models that can automatically generalize to a wide range of cases and problem configurations. This strong generalizability also significantly reduces the need for retraining when encountering new configurations. In other words, the proposed model requires a single training, and provides optimal decisions within a second.

\textit{Randomized training of the proposed model significantly improves model adaptability:}  
    The results showcase that training on datasets with a random number of manufacturing sites, such as \(LRP15M40\) and \(LRP20M50\), significantly enhances the model’s generalizability across diverse configurations. This strategy ensures that the model performs robustly even in dynamic and uncertain operational environments, making it a valuable tool for lessors managing complex, multi-site manufacturing systems.

\section{Summary}
Leased manufacturing systems bring about significant operational and maintenance challenges due to their dynamic and variable nature. The integration of sensor-driven maintenance introduces an additional layer of logistical and computational challenges, as lessors must manage multiple constraints, evolving conditions, and the overarching goal of minimizing costs while maximizing machine reliability and operational efficiency. Traditional optimization approaches, such as MIP, struggle to handle these complexities in real-time due to their computational burden and reliance on expert intervention for problem reformulation.

In this work, we propose an innovative attention-based model, \textit{AttenMfg}, to address these challenges. The model leverages multi-head attention mechanisms to optimize maintenance and operations scheduling in leased manufacturing systems. It is highly flexible, generalizable, and capable of incorporating diverse constraints, such as machine availability, team routing, and operational penalties, ensuring practical feasibility. Our experiments demonstrate that the proposed model outperforms traditional optimization methods in terms of computational efficiency and solution quality, solving problems in seconds where MIP-based approaches fail to converge within hours.

Additionally, the model exhibits strong generalizability, as shown in our experiments where training on larger or more variable datasets improved performance across different scenarios. This capability allows lessors to simulate complex or future settings, train the model only once on these scenarios, and use it to solve a wide range of potential operational problems efficiently and effectively. By doing so, the model not only addresses current challenges but also provides a forward-looking solution for long-term operational planning. 

The proposed attention-based model represents a significant advancement in maintenance and operations scheduling for leased manufacturing systems. It is robust, adaptable, and efficient, making it a practical tool for lessors to manage the complexities of modern manufacturing environments. Its ability to adapt to changing parameters, constraints, and operational settings highlights its potential as a transformative solution for optimizing maintenance and operations in dynamic industrial systems.
\textbf{}







\bibliographystyle{plainnat}
\bibliography{ref}

\begin{thebibliography}{46}
\providecommand{\natexlab}[1]{#1}
\providecommand{\url}[1]{\texttt{#1}}
\expandafter\ifx\csname urlstyle\endcsname\relax
  \providecommand{\doi}[1]{doi: #1}\else
  \providecommand{\doi}{doi: \begingroup \urlstyle{rm}\Url}\fi

\bibitem[Alvarez et~al.(2017)Alvarez, Louveaux, and Wehenkel]{alvarez2017machine}
Alejandro~Marcos Alvarez, Quentin Louveaux, and Louis Wehenkel.
\newblock A machine learning-based approximation of strong branching.
\newblock \emph{INFORMS Journal on Computing}, 29\penalty0 (1):\penalty0 185--195, 2017.

\bibitem[Azadeh et~al.(2015)Azadeh, Farahani, Kalantari, and Zarrin]{azadeh2015solving}
A~Azadeh, M~Hosseinabadi Farahani, SS~Kalantari, and M~Zarrin.
\newblock Solving a multi-objective open shop problem for multi-processors under preventive maintenance.
\newblock \emph{The International Journal of Advanced Manufacturing Technology}, 78:\penalty0 707--722, 2015.

\bibitem[Baltean-Lugojan et~al.(2018)Baltean-Lugojan, Bonami, Misener, and Tramontani]{baltean2018selecting}
Radu Baltean-Lugojan, Pierre Bonami, Ruth Misener, and Andrea Tramontani.
\newblock Selecting cutting planes for quadratic semidefinite outer-approximation via trained neural networks.
\newblock \emph{URL: http://www.optimization-online.org/DB\_HTML/2018/11/6943.html}, 2018.

\bibitem[Bedford et~al.(2011)Bedford, Dewan, Meilijson, and Zitrou]{bedford2011signal}
Tim Bedford, Isha Dewan, Isaac Meilijson, and Athena Zitrou.
\newblock The signal model: A model for competing risks of opportunistic maintenance.
\newblock \emph{European Journal of Operational Research}, 214\penalty0 (3):\penalty0 665--673, 2011.

\bibitem[Camci(2015)]{camci2015maintenance}
Fatih Camci.
\newblock Maintenance scheduling of geographically distributed assets with prognostics information.
\newblock \emph{European Journal of Operational Research}, 245\penalty0 (2):\penalty0 506--516, 2015.

\bibitem[Chang and Lo(2011)]{chang2011joint}
Wen~Liang Chang and Hui-Chiung Lo.
\newblock Joint determination of lease period and preventive maintenance policy for leased equipment with residual value.
\newblock \emph{Computers \& Industrial Engineering}, 61\penalty0 (3):\penalty0 489--496, 2011.

\bibitem[Chen(2009)]{chen2009minimizing}
Wen-Jinn Chen.
\newblock Minimizing number of tardy jobs on a single machine subject to periodic maintenance.
\newblock \emph{Omega}, 37\penalty0 (3):\penalty0 591--599, 2009.

\bibitem[Chen et~al.(2020)Chen, An, Zhang, and Li]{chen2020approximate}
Xiaohui Chen, Youjun An, Zhiyao Zhang, and Yinghe Li.
\newblock An approximate nondominated sorting genetic algorithm to integrate optimization of production scheduling and accurate maintenance based on reliability intervals.
\newblock \emph{Journal of Manufacturing Systems}, 54:\penalty0 227--241, 2020.

\bibitem[Coleman et~al.(2017)Coleman, Damodaran, Chandramouli, and Deuel]{deloitte_maintenance_2017}
Chris Coleman, Satish Damodaran, Mahesh Chandramouli, and Ed~Deuel.
\newblock Making maintenance smarter: Predictive maintenance and the digital supply network.
\newblock Technical report, Deloitte, 2017.
\newblock URL \url{https://www2.deloitte.com/content/dam/insights/us/articles/3828- Making-maintenance-smarter/DUP-Making-maintenance-smarter.pdf}.

\bibitem[Company(2023)]{mckinsey_indirect_2023}
McKinsey~\& Company.
\newblock Indirect manufacturing costs: An overlooked source for clear savings.
\newblock Technical report, McKinsey \& Company, 2023.
\newblock URL \url{https://www.mckinsey.com/capabilities/operations/our-insights/indi rect-manufacturing-costs-an-overlooked-source-for-clear-savings}.

\bibitem[Feng et~al.(2023)Feng, Zhang, Sun, Guo, Fan, Ren, Song, and Wang]{feng2023multi}
Qiang Feng, Yue Zhang, Bo~Sun, Xing Guo, Donming Fan, Yi~Ren, Yanjie Song, and Zili Wang.
\newblock Multi-level predictive maintenance of smart manufacturing systems driven by digital twin: a matheuristics approach.
\newblock \emph{Journal of Manufacturing Systems}, 68:\penalty0 443--454, 2023.

\bibitem[Gebraeel et~al.(2005)Gebraeel, Lawley, Li, and Ryan]{gebraeel2005residual}
Nagi~Z Gebraeel, Mark~A Lawley, Rong Li, and Jennifer~K Ryan.
\newblock Residual-life distributions from component degradation signals: A bayesian approach.
\newblock \emph{IiE Transactions}, 37\penalty0 (6):\penalty0 543--557, 2005.

\bibitem[Geurtsen et~al.(2023)Geurtsen, Adan, and Ak{\c{c}}ay]{geurtsen2023integrated}
Michael Geurtsen, Jelle Adan, and Alp Ak{\c{c}}ay.
\newblock Integrated maintenance and production scheduling for unrelated parallel machines with setup times.
\newblock \emph{Flexible Services and Manufacturing Journal}, pages 1--34, 2023.

\bibitem[Gunarathna et~al.(2022)Gunarathna, Borovica-Gajic, Karunasekara, and Tanin]{gunarathna2022solving}
Udesh Gunarathna, Renata Borovica-Gajic, Shanika Karunasekara, and Egemen Tanin.
\newblock Solving dynamic graph problems with multi-attention deep reinforcement learning.
\newblock \emph{arXiv preprint arXiv:2201.04895}, 2022.

\bibitem[Jia and Zhang(2020)]{jia2020joint}
Chuanzhou Jia and Chi Zhang.
\newblock Joint optimization of maintenance planning and workforce routing for a geographically distributed networked infrastructure.
\newblock \emph{Iise Transactions}, 52\penalty0 (7):\penalty0 732--750, 2020.

\bibitem[Jin et~al.(2023)Jin, Mi, Song, and Li]{jin2023deep}
Xin Jin, Nan Mi, Wen Song, and Qiqiang Li.
\newblock Deep reinforcement learning for dynamic twin automated stacking cranes scheduling problem.
\newblock \emph{Electronics}, 12\penalty0 (15):\penalty0 3288, 2023.

\bibitem[Jung et~al.(2017)Jung, Zhang, and Winslett]{jung2017vibration}
Deokwoo Jung, Zhenjie Zhang, and Marianne Winslett.
\newblock Vibration analysis for iot enabled predictive maintenance.
\newblock In \emph{2017 ieee 33rd international conference on data engineering (icde)}, pages 1271--1282. IEEE, 2017.

\bibitem[Karakaya et~al.(2024)Karakaya, Yildirim, Gebraeel, and Xia]{karakaya2024sensor}
{\c{S}}akir Karakaya, Murat Yildirim, Nagi Gebraeel, and Tangbin Xia.
\newblock A sensor-driven operations and maintenance planning approach for large-scale leased manufacturing systems.
\newblock \emph{International Journal of Production Research}, pages 1--18, 2024.

\bibitem[Kazemian et~al.(2024)Kazemian, Yildirim, and Ramanan]{kazemian2024attention}
Iman Kazemian, Murat Yildirim, and Paritosh Ramanan.
\newblock Attention is all you need to optimize wind farm operations and maintenance.
\newblock \emph{arXiv preprint arXiv:2410.24052}, 2024.

\bibitem[Koochaki et~al.(2012)Koochaki, Bokhorst, Wortmann, and Klingenberg]{koochaki2012condition}
Javid Koochaki, Jos~AC Bokhorst, Hans Wortmann, and Warse Klingenberg.
\newblock Condition based maintenance in the context of opportunistic maintenance.
\newblock \emph{International Journal of Production Research}, 50\penalty0 (23):\penalty0 6918--6929, 2012.

\bibitem[Kool et~al.(2018)Kool, Van~Hoof, and Welling]{kool2018attention}
Wouter Kool, Herke Van~Hoof, and Max Welling.
\newblock Attention, learn to solve routing problems!
\newblock \emph{arXiv preprint arXiv:1803.08475}, 2018.

\bibitem[Liu et~al.(2022)Liu, Zhang, Tang, and Yao]{liu2022good}
Shengcai Liu, Yu~Zhang, Ke~Tang, and Xin Yao.
\newblock How good is neural combinatorial optimization?
\newblock \emph{arXiv preprint arXiv:2209.10913}, 2022.

\bibitem[Liu et~al.(2024)Liu, Liu, and Yang]{liu2024optimal}
Yanping Liu, Biyu Liu, and Haidong Yang.
\newblock Optimal pricing and financing strategies for leased equipment considering maintenance and lessees’ options.
\newblock \emph{International Journal of Production Economics}, 269:\penalty0 109157, 2024.

\bibitem[Liu et~al.(2025)Liu, Liu, Yang, and Luo]{liu2025optimal}
Yanping Liu, Biyu Liu, Haidong Yang, and Kai Luo.
\newblock Optimal production and maintenance strategies for manufacturing/remanufacturing leasing system considering uncertain quality and carbon emission.
\newblock \emph{International Journal of Production Economics}, 280:\penalty0 109489, 2025.

\bibitem[Lyu et~al.(2024)Lyu, Hong, Zhang, and Wang]{lyu2024product}
Zerong Lyu, Zhaofu Hong, Yunrong Zhang, and Chang Wang.
\newblock Product co-design with consumer participation in a service-oriented manufacturing system.
\newblock \emph{International Journal of Production Research}, 62\penalty0 (20):\penalty0 7504--7524, 2024.

\bibitem[Naderi et~al.(2009)Naderi, Zandieh, and Fatemi~Ghomi]{naderi2009scheduling}
B~Naderi, M~Zandieh, and SMT Fatemi~Ghomi.
\newblock Scheduling sequence-dependent setup time job shops with preventive maintenance.
\newblock \emph{The International Journal of Advanced Manufacturing Technology}, 43:\penalty0 170--181, 2009.

\bibitem[Polotski et~al.(2019)Polotski, Kenne, and Gharbi]{polotski2019joint}
V~Polotski, J-P Kenne, and A~Gharbi.
\newblock Joint production and maintenance optimization in flexible hybrid manufacturing--remanufacturing systems under age-dependent deterioration.
\newblock \emph{International Journal of Production Economics}, 216:\penalty0 239--254, 2019.

\bibitem[Rivera-G{\'o}mez et~al.(2016)Rivera-G{\'o}mez, Gharbi, Kenn{\'e}, Montano-Arango, and Hernandez-Gress]{rivera2016production}
H{\'e}ctor Rivera-G{\'o}mez, Ali Gharbi, Jean-Pierre Kenn{\'e}, Oscar Montano-Arango, and Eva~Selene Hernandez-Gress.
\newblock Production control problem integrating overhaul and subcontracting strategies for a quality deteriorating manufacturing system.
\newblock \emph{International Journal of Production Economics}, 171:\penalty0 134--150, 2016.

\bibitem[Rivera-G{\'o}mez et~al.(2018)Rivera-G{\'o}mez, Gharbi, Kenn{\'e}, Monta{\~n}o-Arango, and Hern{\'a}ndez-Gress]{rivera2018subcontracting}
H{\'e}ctor Rivera-G{\'o}mez, Ali Gharbi, Jean-Pierre Kenn{\'e}, Oscar Monta{\~n}o-Arango, and Eva~Selene Hern{\'a}ndez-Gress.
\newblock Subcontracting strategies with production and maintenance policies for a manufacturing system subject to progressive deterioration.
\newblock \emph{International Journal of Production Economics}, 200:\penalty0 103--118, 2018.

\bibitem[Schutz and Rezg(2013)]{schutz2013maintenance}
J{\'e}r{\'e}mie Schutz and Nidhal Rezg.
\newblock Maintenance strategy for leased equipment.
\newblock \emph{Computers \& Industrial Engineering}, 66\penalty0 (3):\penalty0 593--600, 2013.

\bibitem[Sharafali et~al.(2019)Sharafali, Tarakci, Kulkarni, and Hameed]{sharafali2019optimal}
Moosa Sharafali, Hakan Tarakci, Shailesh Kulkarni, and Raja Abdul Razack~Shahul Hameed.
\newblock Optimal delivery due date for a supplier with an unreliable machine under outsourced maintenance.
\newblock \emph{International Journal of Production Economics}, 208:\penalty0 53--68, 2019.

\bibitem[Si et~al.(2022{\natexlab{a}})Si, Xia, Gebraeel, Wang, Pan, and Xi]{si2022reliability}
Guojin Si, Tangbin Xia, Nagi Gebraeel, Dong Wang, Ershun Pan, and Lifeng Xi.
\newblock A reliability-and-cost-based framework to optimize maintenance planning and diverse-skilled technician routing for geographically distributed systems.
\newblock \emph{Reliability Engineering \& System Safety}, 226:\penalty0 108652, 2022{\natexlab{a}}.

\bibitem[Si et~al.(2022{\natexlab{b}})Si, Xia, Pan, and Xi]{si2022service}
Guojin Si, Tangbin Xia, Ershun Pan, and Lifeng Xi.
\newblock Service-oriented global optimization integrating maintenance grouping and technician routing for multi-location multi-unit production systems.
\newblock \emph{IISE transactions}, 54\penalty0 (9):\penalty0 894--907, 2022{\natexlab{b}}.

\bibitem[Taghipour and Banjevic(2012)]{taghipour2012optimum}
Sharareh Taghipour and Dragan Banjevic.
\newblock Optimum inspection interval for a system under periodic and opportunistic inspections.
\newblock \emph{Iie Transactions}, 44\penalty0 (11):\penalty0 932--948, 2012.

\bibitem[Vaswani et~al.(2017)Vaswani, Shazeer, Parmar, Uszkoreit, Jones, Gomez, Kaiser, and Polosukhin]{vaswani2017attention}
Ashish Vaswani, Noam Shazeer, Niki Parmar, Jakob Uszkoreit, Llion Jones, Aidan~N Gomez, {\L}ukasz Kaiser, and Illia Polosukhin.
\newblock Attention is all you need.
\newblock \emph{Advances in Neural Information Processing Systems}, 30, 2017.

\bibitem[Vinyals et~al.(2015)Vinyals, Fortunato, and Jaitly]{vinyals2015pointer}
Oriol Vinyals, Meire Fortunato, and Navdeep Jaitly.
\newblock Pointer networks.
\newblock \emph{Advances in Neural Information Processing Systems}, 28, 2015.

\bibitem[Wan et~al.(2024)Wan, Xu, Chen, and Zhou]{wan2024deep}
Pengfu Wan, Gangyan Xu, Jiawei Chen, and Yaoming Zhou.
\newblock Deep reinforcement learning enabled multi-uav scheduling for disaster data collection with time-varying value.
\newblock \emph{IEEE Transactions on Intelligent Transportation Systems}, 2024.

\bibitem[Wang and Yu(2010)]{wang2010effective}
Shijin Wang and Jianbo Yu.
\newblock An effective heuristic for flexible job-shop scheduling problem with maintenance activities.
\newblock \emph{Computers \& Industrial Engineering}, 59\penalty0 (3):\penalty0 436--447, 2010.

\bibitem[Xia et~al.(2012)Xia, Xi, Zhou, and Lee]{xia2012dynamic}
Tangbin Xia, Lifeng Xi, Xiaojun Zhou, and Jay Lee.
\newblock Dynamic maintenance decision-making for series--parallel manufacturing system based on mam--mtw methodology.
\newblock \emph{European Journal of Operational Research}, 221\penalty0 (1):\penalty0 231--240, 2012.

\bibitem[Xia et~al.(2015)Xia, Jin, Xi, and Ni]{xia2015production}
Tangbin Xia, Xiaoning Jin, Lifeng Xi, and Jun Ni.
\newblock Production-driven opportunistic maintenance for batch production based on mam--apb scheduling.
\newblock \emph{European Journal of Operational Research}, 240\penalty0 (3):\penalty0 781--790, 2015.

\bibitem[Yang et~al.(2008)Yang, Djurdjanovic, and Ni]{yang2008maintenance}
Zimin Yang, Dragan Djurdjanovic, and Jun Ni.
\newblock Maintenance scheduling in manufacturing systems based on predicted machine degradation.
\newblock \emph{Journal of intelligent manufacturing}, 19:\penalty0 87--98, 2008.

\bibitem[Yeh et~al.(2009)Yeh, Kao, and Chang]{yeh2009optimal}
Ruey~Huei Yeh, Kow-Chin Kao, and Wen~Liang Chang.
\newblock Optimal preventive maintenance policy for leased equipment using failure rate reduction.
\newblock \emph{Computers \& Industrial Engineering}, 57\penalty0 (1):\penalty0 304--309, 2009.

\bibitem[Yildirim et~al.(2016{\natexlab{a}})Yildirim, Sun, and Gebraeel]{yildirim2016sensor}
Murat Yildirim, Xu~Andy Sun, and Nagi~Z Gebraeel.
\newblock Sensor-driven condition-based generator maintenance scheduling—part i: Maintenance problem.
\newblock \emph{IEEE Transactions on Power Systems}, 31\penalty0 (6):\penalty0 4253--4262, 2016{\natexlab{a}}.

\bibitem[Yildirim et~al.(2016{\natexlab{b}})Yildirim, Sun, and Gebraeel]{yildirim2016sensor1}
Murat Yildirim, Xu~Andy Sun, and Nagi~Z Gebraeel.
\newblock Sensor-driven condition-based generator maintenance scheduling—part {I}: Maintenance problem.
\newblock \emph{IEEE Transactions on Power Systems}, 31\penalty0 (6):\penalty0 4253--4262, 2016{\natexlab{b}}.

\bibitem[Zhou et~al.(2009)Zhou, Xi, and Lee]{zhou2009opportunistic}
Xiaojun Zhou, Lifeng Xi, and Jay Lee.
\newblock Opportunistic preventive maintenance scheduling for a multi-unit series system based on dynamic programming.
\newblock \emph{International Journal of Production Economics}, 118\penalty0 (2):\penalty0 361--366, 2009.

\bibitem[Zhou et~al.(2012)Zhou, Lu, and Xi]{zhou2012preventive}
Xiaojun Zhou, Zhiqiang Lu, and Lifeng Xi.
\newblock Preventive maintenance optimization for a multi-component system under changing job shop schedule.
\newblock \emph{Reliability Engineering \& System Safety}, 101:\penalty0 14--20, 2012.

\end{thebibliography}

\end{document}